\documentclass{aip-cp}

\usepackage[numbers]{natbib}
\usepackage{rotating}
\usepackage{graphicx}

\newcommand{\df}{\displaystyle\frac}
\newcommand{\la}{\lambda}
\newcommand{\al}{\alpha}
\newcommand{\om}{\omega}

\newcommand{\be}{\beta}

\newtheorem{proposition}{Proposition}

\begin{document}

\title{Stability analysis of a hypothalamic-pituitary-adrenal axis model with inclusion of glucocorticoid receptor and memory}

\author[aff1,aff2]{Eva Kaslik}
\eaddress{ekaslik@gmail.com}
\author[aff3]{Dan Bogdan Navolan}
\eaddress{navolan@yahoo.com}
\author[aff1,aff4]{Mihaela Neam\c{t}u\corref{cor1}}
\eaddress{mihaela.neamtu@e-uvt.ro}

\affil[aff1]{West University of Timisoara, Bd. V. Parvan nr. 4, 300223, Romania}
\affil[aff2]{Institute e-Austria Timisoara, Bd. V. Parvan nr. 4, cam 045B, 300223, Romania}
\affil[aff3]{University of Medicine and Pharmacy "Victor Babes", 2 P-ta Eftimie Murgu, Timisoara, Romania}
\affil[aff4]{Politehnica University of Bucharest, 313 Splaiul Independentei, 060042 Bucharest, Romania}
\corresp[cor1]{Corresponding author: mihaela.neamtu@e-uvt.ro}

\maketitle

\begin{abstract}
This paper analyzes a four-dimensional model of the hypothalamic-pituitary-adrenal (HPA) axis that includes the influence of the glucocorticoid receptor in the pituitary. Due to the spatial separation between the hypothalamus, pituitary and adrenal glands, distributed time delays are introduced in the mathematical model. The existence of the positive equilibrium point is proved and a local stability and bifurcation analysis is provided, considering several types of delay kernels.  The fractional-order model with discrete time delays is also taken into account. Numerical simulations are provided to illustrate the effectiveness of the theoretical findings.
\end{abstract}

\section{INTRODUCTION}

One of the most important systems in stress response is hypothalamus-pituitary-adrenal (HPA) axis, because the released hormones lead to energy directed to the organisms. The HPA axis is organized into three distinct regions: the hypothalamus, pituitary gland and adrenal gland. It is a complex set of direct influences and feedback interactions among the three endocrine glands. These glands work together by producing and secreting, or responding to common hormones including corticotropin-releasing hormone (CRH), adenocorticotropin  (ACTH), and cortisol (CORT) \cite{Kyrylov_2005}. The dynamics of the HPA axis is very important, because cortisol overproduction leads to Cushing’s disease;
cortisol underproduction generates Addison’s disease.

During the past decade, the mathematical modelling of the HPA axis has been intensively studied \cite{Kyrylov_2005,Andersen_2013,Bairagi_2008,Conrad_2009,Gudmand_2014,Jelic_2005,Lenbury_2005,Markovic_2011,Pornsawad_2013,Savic_2005,Savic_2006,Vinther_2011}.
In 2005, Savic and Jelic \cite{Savic_2005} presented one of the mathematical models of HPA axis described by a system with three nonlinear differential equations with one delay.  Three variables have been considered: the concentration of $CRH$, the concentration of $ACTH$ and the concentration of free cortisol $CORT$. According to the results obtained by stability analysis, the system does not oscillate at all in the neighborhood of the equilibrium point. In 2006,  Savic et al. \cite{Savic_2006} considered the same model, but four time delays have been introduced. Once more, their results showed that the system describing the HPA axis does not oscillate in a neighborhood of the equilibrium points.

In 2011, Vinther et al. \cite{Vinther_2011} introduced feedback functions as generalized Hill functions, based on the receptor dynamics. They investigated the possibility of time delays as being capable of producing oscillations in accordance with the ultradian rhythm. By numercial simulations, oscillating solutions in a neighborhood of the equilibrium point have been observed.

In 2015, Kaslik and Neamtu \cite{EvaMiha_MMB} proved the occurrence of oscillating solutions in a neighborhood of the unique equilibrium point, by performing a Hopf bifurcation analysis, considering several distributed time delays and fractional derivatives in the minimal three-dimensional mathematical model previously analyzed in \cite{Vinther_2011}.

On the other hand, Gupta et al. \cite{Gupta_2007} and Sriram et al. \cite{Sriram_2012} incorporated expression of the glucocorticoid receptor (GR) in the mathematical model of the HPA axis and demonstrated that repeated stress and the level of GR concentration leads to a disorder; however, it is important to note that they did not consider time delays.

With the aim of improving the mathematical model of the HPA axis described in \cite{Gupta_2007}, we include distributed time delays and  frac\-tional-order derivatives. On one hand, distributed time delays represent the situation where the delays occur in certain ranges of values with some associated probability distributions, taking into account the whole past history of the variables. In many real world applications, distributed time delays are more realistic and more accurate than discrete time delays \cite{Cushing_2013}. Distributed delay models appear in a wide range of applications such as, population biology \cite{Faria_2008,Ruan_1996}, hematopoiesis \cite{Adimy_2003,Adimy_2005,Adimy_2006,Ozbay_2008}, neural networks \cite{Jessop_2010}. On the other hand, the main benefit of fractional-order models in comparison with classical integer-order models is that fractional derivatives provide a good tool for the description of memory and hereditary properties of various processes \cite{Kilbas,Lak,Podlubny}.

The paper is structured as follows. Section 2 provides the mathematical model of the HPA axis, with the influence of the GR concentration, where we introduce distributed time delays to account for the time needed by the hormones to travel from source to destination, as well as fractional-order derivatives to account for the memory properties of the system. In Section 3, a local stability analysis of the system with distributed delays is provided. Numerical simulations are carried out and discussed in Section 4, followed by concluding remarks in Section 5.

\section{THE MATHEMATICAL MODEL}
The mechanism of the HPA axis can be shortly described as follows. Paraventricular nuclei (PVN) of hypothalamus generate corticotropin releasing hormone (CRH) which induces the adenocorticotropin (ACTH) production in the pituitary. Then, ACTH stimulates the adrenal cortex to produce glucocorticoids, which in turn suppress the production of both CRH and ACTH.

\begin{figure}[htbp]
\centering
\includegraphics*[width=0.4\textwidth]{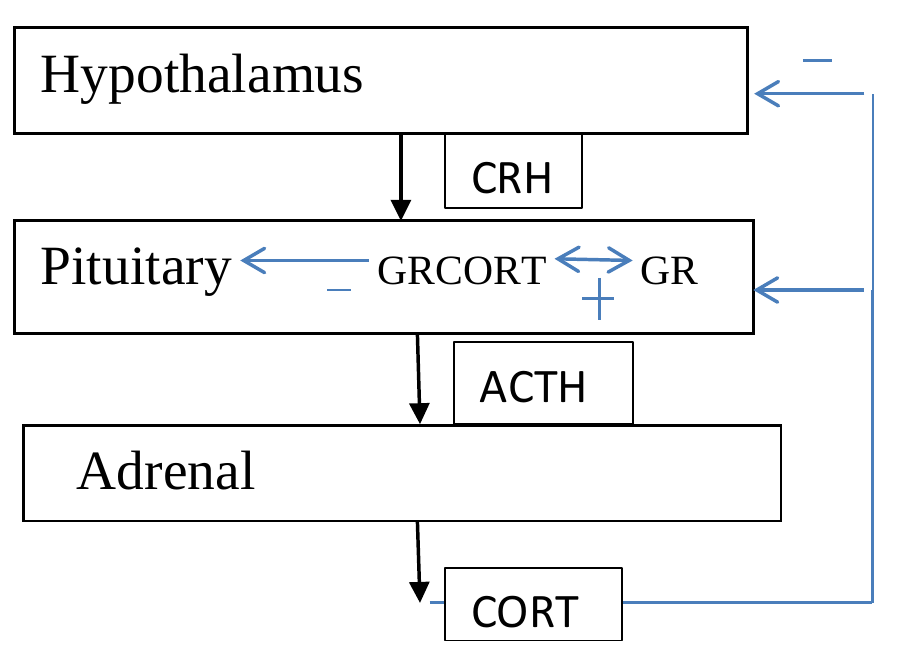}
\caption{A simple schematic representation of the HPA axis with negative feedback.}
\label{fig:0}
\end{figure}

In formulating the mathematical model which describes the variation in time of the concentrations of the three hormones CRH, ACTH and CORT, as well as the glucocorticoid receptors GR, the following sequence of typical events is considered, according to the schema presented in Figure \ref{fig:0}.

CRH is secreted from the hypothalamus and released into the portal blood vessel of the hypophyseal stalk, and then transported to the anterior pituitary where it stimulates the secretion of ACTH. In our mathematical model, the evolution of the CRH hormone concentration in the hypothalamus is described by:
$$\dot CRH(t)=\df{a_1}{a_2+\int_{-\infty}^th_1(t-s)CORT(s)ds}-a_3CRH(t),$$
where $a_1$ models a circadian production and stress term, $a_2$ is the inhibition constant, $a_3$ is the degradation rate of cortisol.

In the pituitary, cortisol enters the cells and binds the glucocorticoid receptor in the cytoplasm, causing the receptor to dimerize. This leads to the complex $GR\times CORT$ that inhibits the production of $ACTH$, modeled by:
$$\dot ACTH(t)=\df{b_1CRH(t)}{b_2+GR(t)\int_{-\infty}^th_1(t-s)CORT(s)ds}-b_3ACTH(t),$$
where $b_1$ stands for a production term, $b_2$ is the inhibition constant, $b_3$ is the degradation rate of $ACTH$, and
$$\dot GR(t)=\df{c_1\left(GR(t)\int_{-\infty}^th_1(t-s)CORT(s)ds\right)^2}{c_2+\left (GR(t)\int_{-\infty}^th_1(t-s)CORT(s)ds\right)^2}+c_3-c_4GR(t),$$
where $c_1$ is the dimerization constant, $c_2$ the binding affinity constant, $c_3$ the production constant and $c_4$ the degradation term for pituitary $GR$ production.

As for the evolution of cortisol concentration produced by the adrenal, the following equation is considered:
$$\dot CORT(t)=d_1\int_{-\infty}^th_2(t-s)ACTH(s)ds-d_2CORT(t),$$
where $d_1$ is a production coefficient and $d_2$ is the degradation rate of $CORT$.

Therefore, the four-dimensional mathematical model with distributed time delay is:
\begin{equation}\label{sysHPA4}
\begin{array}{l}
\dot CRH(t)=\df{a_1}{a_2+\int_{-\infty}^th_1(t-s)CORT(s)ds}-a_3CRH(t),\\
\dot ACTH(t)=\df{b_1CRH(t)}{b_2+GR(t)\int_{-\infty}^th_1(t-s)CORT(s)ds}-b_3ACTH(t),\\
\dot GR(t)=\df{c_1\left(GR(t)\int_{-\infty}^th_1(t-s)CORT(s)ds\right)^2}{c_2+\left (GR(t)\int_{-\infty}^th_1(t-s)CORT(s)ds\right)^2}+c_3-c_4GR(t),\\
\dot CORT(t)=d_1\int_{-\infty}^th_2(t-s)ACTH(s)ds-d_2CORT(t),\\
\end{array}
\end{equation}
with the initial conditions $CRH(t)=\phi_1(t)$, $ACTH(t)=\phi_2(t)$, $GR(t)=\phi_3(t)$, $CORT(t)=\phi_4(t)$, $t\in(-\infty,0]$, where $\varphi_i$  are bounded continuous functions defined on $(-\infty,0]$, with values in $[0,\infty)$.

In system (\ref{sysHPA4}), the delay kernels $h_1,h_2:[0,\infty)\to[0,\infty)$ are probability density functions representing the probability that a particular time delay occurs. They are assumed to be bounded, piecewise continuous and satisfy
\begin{equation}\label{delay.kernel.properties}
\int_0^{\infty}h(s)ds=1.
\end{equation}
The average delay of a delay kernel $h(t)$ is given by
$$\tau=\int_0^{\infty}sh(s)ds<\infty.$$
Two important classes of delay kernels often used in the literature are worth mentioning:
\begin{itemize}
\item Dirac kernels: $h(s)=\delta(s-\tau)$, where $\tau\geq 0$. In this particular case, the distributed delay is reduced to a discrete time delay:
$$\int_{-\infty}^t x(s)h(t-s)ds=\int_0^\infty x(t-s)\delta(s-\tau)ds=x(t-\tau).$$
\item Gamma kernels: $h(s)=\df{ s^{p-1}e^{-s/\beta}}{\beta^p\Gamma(p)}$, where $p,\beta>0$. The average delay of a Gamma kernel is $\tau=p\beta$.
\end{itemize}

On the other hand, we also consider the four-dimensional mathematical model with fractional-order derivatives of Caputo type and discrete time delays $\tau_1,\tau_2>0$:
\begin{equation}\label{sysHPA4.frac}
\begin{array}{l}
^cD^q CRH(t)=\df{a_1}{a_2+CORT(t-\tau_1)}-a_3CRH(t),\\
^cD^q ACTH(t)=\df{b_1CRH(t)}{b_2+GR(t)CORT(t-\tau_1)}-b_3ACTH(t),\\
^cD^q GR(t)=\df{c_1\left(GR(t)CORT(t-\tau_1)\right)^2}{c_2+\left (GR(t)CORT(t-\tau_1)\right)^2}+c_3-c_4GR(t),\\
^cD^q CORT(t)=d_1ACTH(t-\tau_2)-d_2CORT(t),\\
\end{array}
\end{equation}
with the initial conditions $CRH(t)=\phi_1(t)$, $ACTH(t)=\phi_2(t)$, $GR(t)=\phi_3(t)$, $CORT(t)=\phi_4(t)$, $t\in(-\max(\tau_1,\tau_2),0]$, where $\varphi_i$  are bounded continuous functions defined on $(-\max(\tau_1,\tau_2),0]$, with values in $[0,\infty)$.

\section{STABILITY OF THE EQUILIBRIUM POINT}

In what follows we use the notations: $CRH(t)=x_1(t)$, $ACTH(t)=x_2(t)$, $GR(t)=x_3(t)$, $CORT(t)=x_4(t)$. The coordinates of the equilibrium point of systems (\ref{sysHPA4}) and (\ref{sysHPA4.frac}) are the positive solution of the following algebraic system:
$$
\left\{
\begin{array}{l}
a_1-a_2a_3x_1-a_3x_1x_4=0,\\
b_1x_1-b_2b_3x_2-b_3x_2x_3x_4=0,\\
c_4x_3(x_2x_4)^2-(c_1+c_3)(x_3x_4)^2+c_2c_4x_3-c_2c_3=0,\\
d_1x_2-d_2x_4=0.
\end{array}\right.
$$

\begin{proposition} $ $
\begin{itemize}
\item[(i)] The equation
$$
\begin{array}{l}
F(x):=c_4d_1^2(p_2x^2+p_1x+p_0)^3-(c_1+c_3)d_1^2x^2(q_1x+q_0)(p_2x^2+p_1x+p_0)^2+\\
\qquad\quad+c_2c_4d_2^2x^2(q_1x+q_0)^2(p_2x^2+p_1x+p_0)-c_2c_3d_2^2x^4(q_1x+q_0)^3=0,
\end{array}
$$
has at least one positive root;

\item[(ii)] An equilibrium point $E_0$ of systems (\ref{sysHPA4}) and (\ref{sysHPA4.frac}) has the coordinates $(x_{10}$, $x_{20}$, $x_{30}$, $x_{40})$, where $x_{20}$ is a positive solution of equation $F(x)=0$ and
$$x_{10}=\df{a_1d_2}{d_2a_2a_3+d_1a_3x_{20}}\quad;\quad x_{30}=\df{p_2x_{20}^2+p_1x_{20}+p_0}{x_{20}^2(q_1x_{20}+q_0)}\quad;\quad
x_{40}=\df{d_1x_{20}}{d_4}
$$
\end{itemize}
\end{proposition}

In the case of the distributed time-delay system (\ref{sysHPA4}), considering the transformation:
$$
u_1(t)=x_1(t)-x_{10},\quad u_2(t)=x_2(t)-x_{20},\quad u_3(t)=x_3(t)-x_{30},\quad u_4(t)=x_4(t)-x_{40},
$$
the linearized system at an equilibrium point $E_0$ becomes:
\begin{equation}\label{14}
\begin{array}{l}
\dot u_1(t)=a_{11}u_1(t)+b_{14}\int_{-\infty}^th_1(s)u_4(t-s)ds,\\
\dot u_2(t)=a_{21}u_1(t)+a_{22}u_2(t)+a_{23}u_3(t)+b_{24}\int_{-\infty}^th_1(s)u_4(t-s)ds,\\
\dot u_3(t)=a_{33}u_3(t)+b_{34}\int_{-\infty}^th_1(s)u_4(t-s)ds,\\
\dot u_4(t)=b_{42}\int_{-\infty}^th_2(s)u_2(t-s)ds+a_{44}u_4(t),
\end{array}
\end{equation} where
$a_{11}=-a_3$, $b_{14}=-\df{a_1}{(a_2+x_{40})^2}$, $a_{21}=\df{b_1}{b_2+x_{30}x_{40}}$, $a_{22}=-b_3$, $a_{23}=-\df{b_1x_{10}x_{40}}{(b_2+x_{30}x_{40})^2}$, $b_{24}=-\df{b_1x_{10}x_{30}}{(b_2+x_{30}x_{40})^2}$, $a_{33}=-c_4+\df{2c_1c_2x_{30}x_{40}^2}{(c_2+x_{30}^2x_{40}^2)^2}$,  $b_{34}=\df{2c_1c_2x_{30}^2x_{40}}{(c_2+x_{30}^2x_{40}^2)^2}$, $b_{42}=d_1$, $a_{44}=-d_2$.

The characteristic equation of (\ref{14}) is given by:
\begin{equation}\label{16}
\la^4+r_{3}\la^3+r_{2}\la^2+r_{1}\la+r_{0}+(s_{2}\la^2+s_{1}\la+s_{0})\left(\int_{-\infty}^0h_1(-s)e^{\lambda s}ds\right)\left(\int_{-\infty}^0h_2(-s)e^{\lambda s}ds\right)=0,
\end{equation} where
$r_3=-a_{11}-a_{22}-a_{33}-a_{44}$, $r_2=a_{11}(a_{22}+a_{33}+a_{44})+a_{33}a_{44}+a_{22}(a_{33}+a_{44})$, $r_1=-a_{33}a_{44}(a_{11}+a_{22})-a_{11}a_{22}(a_{33}+a_{44})$, $r_0=a_{11}a_{22}a_{33}a_{44}$, $s_2=-b_{42}b_{24}$, $s_1=-b_{42}(a_{21}b_{14}+a_{23}b_{34}-b_{24}(a_{11}+a_{33}))$, $s_0=-b_{42}(-a_{23}b_{34}a_{11}-a_{21}b_{14}a_{33}+b_{24}a_{11}a_{33})$.

If there are no delays, the characteristic equation (\ref{16}) is given by:
\begin{equation}\label{160}
\la^4+r_3\la^3+(r_2+s_2)\la^2+(r_1+s_1)\la+r_0+s_0=0
\end{equation}
Using the Hurwitz criteria  for equation (\ref{160}) we obtain:

\begin{proposition}
In the non-delayed case, if inequalities
\begin{equation}
r_3>0,\quad r_1+s_1>0,\quad r_0+s_0>0,\quad r_3(r_2+s_2)(r_1+s_1)>(r_1+s_1)^2+r_3^2(r_0+s_0),
\end{equation}
hold, then the equilibrium point $E_0(x_{10}, x_{20}, x_{30}, x_{40})$ of system (\ref{sysHPA4}) is locally asymptotically stable.
\end{proposition}

In the presence of different types of time delays, we analyze three cases, as follows.

\noindent\textbf{Case 1. Dirac kernels:} $h_i(s)=\delta(s-\tau_i)$, $i=1,2$.

In this case, denoting $\tau=\tau_1+\tau_2$, the characteristic equation (\ref{16}) becomes:
\begin{equation}\label{18}
\la^4+r_{3}\la^3+r_{2}\la^2+r_{1}\la+r_{0}+(s_{2}\la^2+s_{1}\la+s_{0})e^{-\lambda \tau}=0.
\end{equation}
If $i\om$, with $\om>0$, is a root of (\ref{18}) then $\om$ satisfies the following equation:
\begin{equation}\label{19}
\om^4-ir_{3}\om^3-r_{2}\om^2+ir_{1}\om+r_{0}+(s_{2}\om^2-is_{1}\om-s_{0})(\cos(\om\tau)-i\sin(\om\tau))=0.
\end{equation}
Separating the real and imaginary parts of (\ref{19}) we obtain the system:
\begin{equation}\label{20}
\left\{
\begin{array}{l}
\om^4-r_{2}\om^2+r_{0}=(-s_{2}\om^2+s_{0})\cos(\om\tau)+s_{1}\om\sin(\om\tau)\\
-r_{3}\om^3+r_{1}\om=(s_{2}\om^2-s_{0})\sin(\om\tau)+s_{1}\om\cos(\om\tau).
\end{array}
\right.
\end{equation}
From this system, it follows that $\om$ satisfies the equation:
\begin{equation}\label{21}
\om^8+m_6\om^6+m_4\om^4+m_2\om^2+m_0=0,
\end{equation} where $m_6=r_{3}^2-2r_{2}$, $m_4=r_{2}^2+2r_{0}-2r_{3}r_{1}-r_{2}^2$, $m_2=r_{1}^2-2r_{2}r_{0}-s_{1}^2+2s_{0}s_{2}$, $m_0=r_{0}^2-s_0^2$.

Denoting $z=\om^2$, eq. (\ref{21}) can be rewritten as:
\begin{equation}\label{23}
z^4+m_1z^3+m_2z^2+m_0z+m_4=0.
\end{equation}
Let $z_0$ denote the smallest positive real root of equation (\ref{23}) and $\omega_0=\sqrt{z_0}$. Then:
\begin{equation}\label{31}
\cos(\om_0\tau)=\df{\Delta_1}{\Delta},
\end{equation} where
\begin{equation}\label{32}
\begin{array}{l}
\Delta_1=s_{2}\om_0^6+(r_{3}s_{1}-s_{2}r_{3}-s_{0})\om_0^4+(r_{0}s_{2}+r_{3}s_{0}-r_{1}s_{1})\om_0^2-r_{0}s_{0},\\
\Delta=s_{2}^2\om_0^4+(s_{1}^2-2s_{0}s_{2})\om_0^2+s_{0}^2.
\end{array}
\end{equation}
 Therefore,
\begin{equation}\label{33}
\tau_j=\df{1}{\om_0}\left[\arccos\left(\df{\Delta_1}{\Delta}\right)+2j\pi\right], j=0,1,2,...
\end{equation}

\begin{proposition} Suppose that $z_0=\om_0^2$, and $\df{dh(z)}{dz}\not =0$, where
$$h(z)=z^4+m_6z^3+m_4z^2+m_2z+m_0.$$ The sign of $Re\left[\df{d\la(\tau)}{d\tau}|_{\tau=\tau_j}\right]$ is the same as that of $\df{dh(z)}{dz}$. Therefore, the following transversality condition holds:
$$Re\left[\df{d\la(\tau)}{d\tau}|_{\tau=\tau_j}\right]\not=0,$$
which corresponds to the occurrence of a Hopf bifurcation in a neighborhood of the equilibrium point $E_0$.
\end{proposition}

\noindent\textbf{Case 2. Mixed kernels:} $h_1(s)=\delta(s-\tau_1)$, $\tau_1>0$, $h_2(s)=a_{20}e^{-a_{20}s}$, $a_{20}>0$.

In this case equation (\ref{16}) becomes:
\begin{equation}\label{181}
\la^5+\al_{4}\la^4+\al_{3}\la^3+\al_2\la^2+\al_{1}\la+\al_{0}+a_{20}(s_{2}\la^2+s_{1}\la+s_{0})e^{-\lambda \tau_1}=0,
\end{equation} where $\al_4=a_{20}+r_3$, $\al_3=a_{20}r_3+r_2$, $\al_2=a_{20}r_2+r_1$, $\al_1=a_{20}r_1+r_0$, $\al_0=a_{20}r_0$ and $s_2$, $s_1$, $s_0$ above.

If $i\om$, with $\om>0$ is a root of (\ref{181}) then $\om$ satisfies the following equation:
\begin{equation}\label{191}
-i\om^5+\al_4\om^4-i\al_3\om^3-\al_{2}\om^2+i\al_{1}\om+\al_{0}+a_{20}(-s_{2}\om^2+is_{1}\om+s_{0})(\cos(\om\tau_1)-i\sin(\om\tau_1))=0.
\end{equation}
Separating the real and imaginary parts of (\ref{191}) we obtain that $\om$ satisfies the equation:
\begin{equation}\label{211}
\om^{10}+q_8\om^8+q_6\om_6+q_4\om^4+q_2\om^2+q_0=0,
\end{equation} where $q_8=\al_1^2+2\al_3$, $q_6=\al_{3}^2-2\al_{1}-2\al_1\al_2$, $q_4=\al_{2}^2+2\al_{0}-2\al_{3}\al_{1}-a_{20}^2s_{2}^2$, $q_2=\al_{1}^2-2\al_{2}\al_{0}-\al_{1}^2+2a_{20}^2s_{0}s_{2}$, $q_0=\al_{0}^2-a_{20}^2s_0^2$.

Considering $z=\om^2$, eq. (\ref{211}) can be rewritten as:
\begin{equation}\label{231}
z^5+q_8z^4+q_4z^2+m_2z+q_0=0.
\end{equation}
Let $z_0$ denote the smallest positive real root of equation (\ref{231}) and $\omega_0=\sqrt{z_0}$. Then:
\begin{equation}\label{311}
\cos(\om_0\tau_1)=\df{A_1}{a_{20}B_1},
\end{equation} where
\begin{equation}\label{321}
\begin{array}{l}
A_1=(\al_1\om_0^4-\al_2\om_0^2+\al_0)(s_2\om_0^2-s_0)+(\om_0^5+\al_3\om_0^3-\al_1\om_0)s_1\om_0,\\
B_1=(s_1\om_0-s_0)^2+(s_1\om_0)^2.
\end{array}
\end{equation}
 Therefore,
\begin{equation}\label{331}
\tau_{1j}=\df{1}{\om_0}\left[arccos\left(\df{A_1}{a_{20}B_1}\right)+2j\pi\right],  j=0,1,2,...
\end{equation}

\begin{proposition} Suppose that $z_0=\om^2_0$,  and $\df{dh(z)}{dz}\not =0$, where
$$h(z)=z^5+q_8z^4+q_4z^2+q_2z+q_0.$$ The sign of $Re\left[\df{d\la(\tau)}{d\tau}|_{\tau=\tau_{1j}}\right]$ is the same as that of $\df{dh(z)}{dz}$. Therefore, the following transversality condition holds:
$$Re\left[\df{d\la(\tau)}{d\tau}|_{\tau=\tau_{1j}}\right]\not=0,$$
which corresponds to the occurrence of a Hopf bifurcation in a neighborhood of the equilibrium point $E_0$.
\end{proposition}

\noindent\textbf{Case 3. Weak gamma kernels:} $h_i(s)=ae^{-as}$, $a>0$, $i=1,2$.

The characteristic equation becomes:
\begin{equation}\label{451}
\la^6+\be_{5}\la^5+\be_{4}\la^4+\be_{3}\la^3+\be_{2}\la^2+\be_{1}\la+\be_0=0,
\end{equation}
where
$\be_5=r_3+2a$, $\be_4=2r_3a+r_2$, $\be_3=2r_2a+r_1+a^2$, $\be_2=r_0+2r_1a+r_2a^2+s_2a^2$, $\be_1=r_1a^2+s_1a^2+2r_0a$, $\be_0=r_0a^2+s_0a^2$.
Denoting
$$
\begin{array}{l}
D_1(a)=\beta_5\beta_4-\beta_3,\\
D_2(a)=\beta_5\beta_4\beta_3+\beta_5\beta_1-\beta_3^2-\beta_2\beta_5,\\
D_3(a)=\beta_5(\beta_5\beta_3\beta_2+\beta_5\beta_4\beta_0+\beta_2\beta_1-\beta_3\beta_0-\beta_5\beta_2^2-\beta_1\beta_4^2)-(\beta_3^2\beta_2+\beta_1^2-\beta_5\beta_2\beta_1-\beta_4\beta_3\beta_1),\\
D_4(a)=\beta_1\beta_3-\beta_0
\left|\begin{array}{cccc}
\beta_5&1&0&0\\
\beta_3&\beta_4&\beta_5&0\\
\beta_1&\beta_2&\beta_3&\beta_5\\
0&\beta_0&\beta_1&\beta_3
\end{array}\right|, \\
D_5(a)=\beta_0D_4(a)\\
\end{array}
$$
we obtain the following:

\begin{proposition}
\begin{itemize}
\item[(i)] If the conditions $D_1(a)>0$, $D_2(a)>0$, $D_3(a)>0$, $D_4(a)>0$, $D_5(a)>0$ hold, for any $a> 0$, the equilibrium point $E_0$  is locally asymptotically stable;

\item[(ii)] If there exists $a_{0}>0$ so that $D_4(a_0)=0$ and $D_1(a_0)>0$, $D_2(a_0)>0$, $D_3(a_0)>0$ and $\df{dD_4(a)}{da}|_{a=a_0}\not =0$, a Hopf bifurcation occurs at $E_0$ as $a$ passes through $a_0$.
\end{itemize}
\end{proposition}

\section{NUMERICAL SIMULATIONS}

For the numerical simulations, we consider the scaled parameter values used in \cite{Gupta_2007}:
$$a_1 = a_2 = 0.1,~~a_3 = 1,~~ b_1 = b_2 = 0.1,~~ b_3 = 10, ~~c_1 = 1,~~c_2 = 0.001,~~ c_3 = 0.05,~~c_4 = 0.9,~~d_1 = d_2 = 1.$$
For these values of the system parameters, three
positive equilibrium states coexist because of the homodimerization
of the GR with cortisol: one equilibrium state with a low GR concentration which is asymptotically stable for any choice of time delays, one equilibrium state with a medium GR concentration which is unstable for any choice of time delays and one equilibrium state with a high GR concentration whose stability depends on the choice of time delays. According to \cite{Gupta_2007}, the low GR concentration represents the normal state, and
high GR concentration reflects a dysregulated HPA axis resulting in persistently low cortisol levels.

In the following, we investigate the stability of the equilibrium point with high GR concentration:\\
$$E_0:~~\bar{x}_{10}=0.66013~(CRH),~~\bar{x}_{20}=0.0514~ (ACTH),~~ \bar{x}_{30}=0.5481 ~(GR),~~\bar{x}_{40}=0.0514~ (CORT).$$

In the case of discrete time-delays: $h_1(s)=\delta(s-\tau_1)$, $h_2(s)=\delta(s-\tau_2)$, with $\tau_1=25$ min, $\tau_2=\tau-\tau_1$ min, we compute the critical value $\tau^\star=32.8043$ min, according to Proposition 3. Figure 2 displays the trajectories of system (\ref{sysHPA4}) for $\tau=\tau_1+\tau_2=31$ min and $\tau=\tau_1+\tau_2=33$ min, respectively. When $\tau=31$ min (below the Hopf bifurcation value), we observe that the trajectories of the system converge to the asymptotically stable equilibrium state $E_0$. On the other hand,  when $\tau=33$ min (above the Hopf bifurcation value), the trajectories of the system converge to the asymptotically stable limit cycle which appears in a neighborhood of $E_0$ due to the supercritical Hopf bifurcation.

\begin{figure}[htbp]
\centering
\begin{tabular}{cc}
\includegraphics[width=0.49\textwidth ]{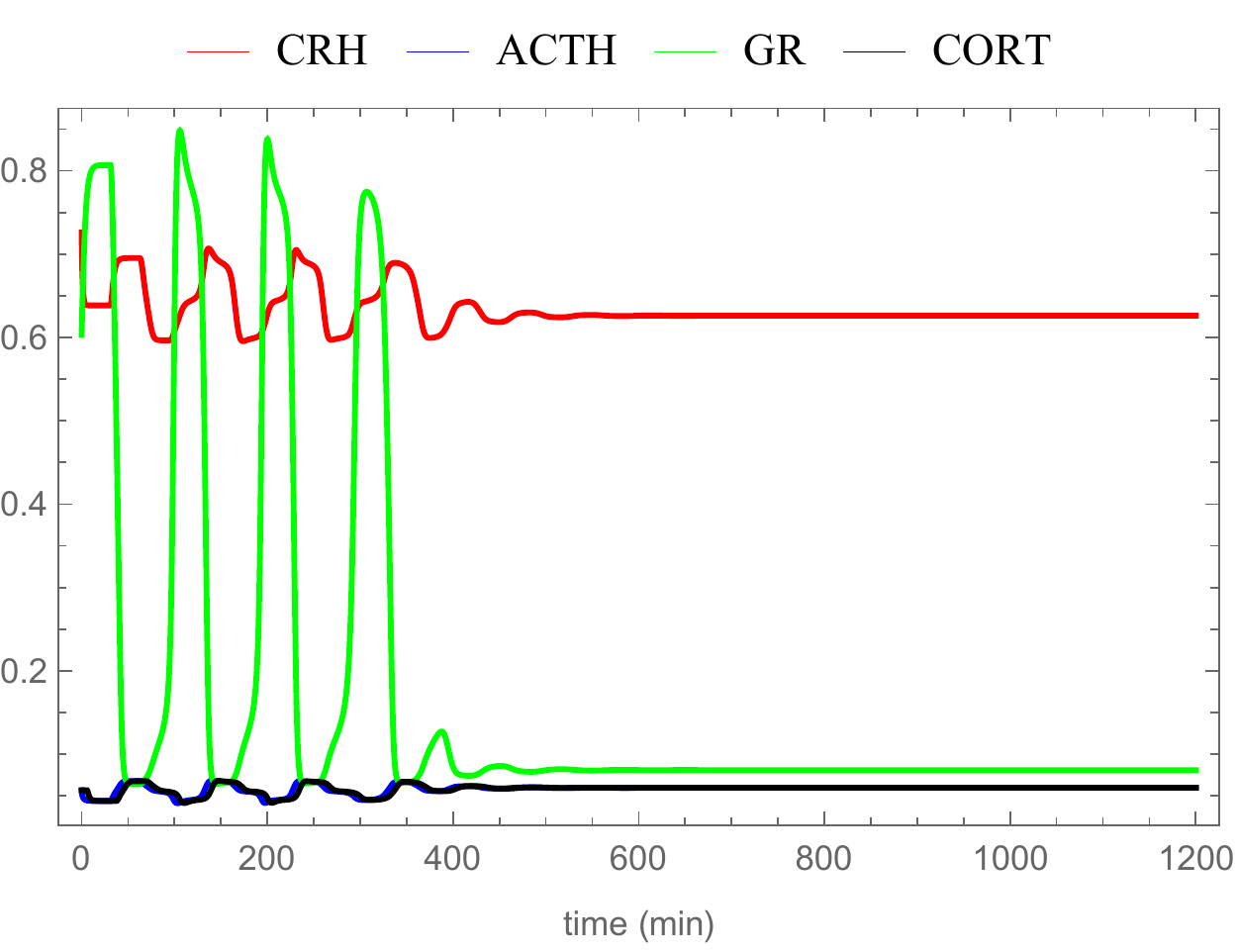} &
\includegraphics[width=0.49\textwidth ]{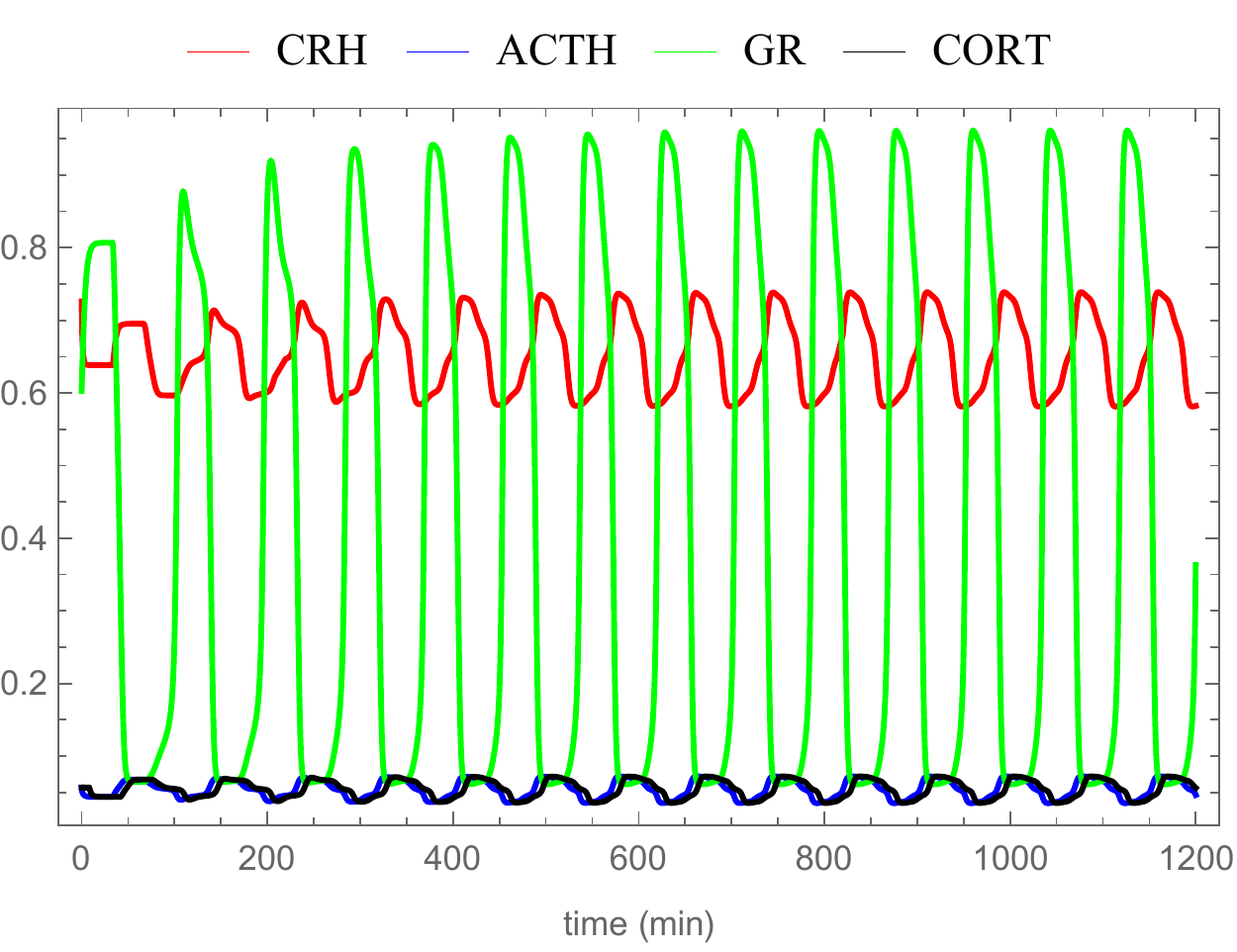}\\
\end{tabular}
\caption{Trajectories of system (\ref{sysHPA4}) with discrete time-delays $\tau_1=25$ min and $\tau_2=6$ min (left) and $\tau_2=8$ min (right), considering an initial condition in a small neighborhood of $E_0$. }
\end{figure}

In the case of weak gamma kernels: $h_1(s)=h_2(s)=ae^{-as}$, with $a=\tau^{-1}$, the equilibrium state $E_0$ is asymptotically stable, for any $a>0$. Figure 3 shows the trajectories of system (\ref{sysHPA4}) for $\tau=50$ min, where the trajectories of the system converge to the asymptotically stable equilibrium state $E_0$.

\begin{figure}[htbp]
\centering
\includegraphics[width=0.49\textwidth ]{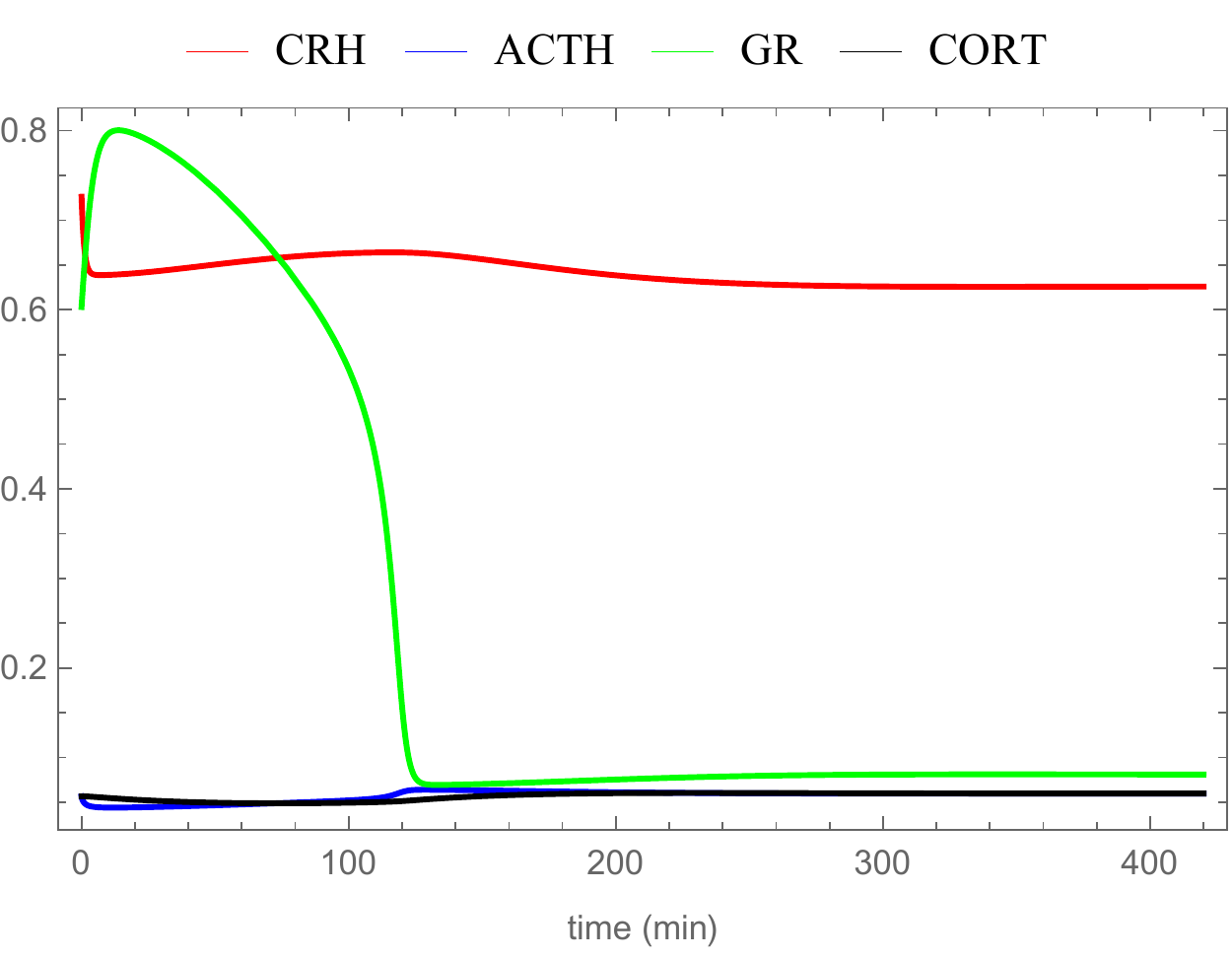}
\caption{Trajectories of system (\ref{sysHPA4}) with weak gamma kernels with $\tau=a^{-1}=50$ min, considering an initial condition in a small neighborhood of $E_0$. }
\end{figure}

For the fractional-order model (\ref{sysHPA4.frac}) with $q=0.8$, we consider the discrete time-delays: $h_1(s)=\delta(s-\tau_1)$, $h_2(s)=\delta(s-\tau_2)$, with $\tau_1=25$ min, $\tau_2=\tau-\tau_1$ min. Figure 4 displays the trajectories of system (\ref{sysHPA4.frac}) for $\tau=\tau_1+\tau_2=39$ min and $\tau=\tau_1+\tau_2=40$ min, respectively. When $\tau=39$ min (below the Hopf bifurcation value), we observe that the trajectories of the system converge to the asymptotically stable equilibrium state $E_0$. On the other hand,  when $\tau=40$ min (above the Hopf bifurcation value), the trajectories of the system converge to the asymptotically stable limit cycle which appears in a neighborhood of $E_0$ due to the supercritical Hopf bifurcation. Numerical simulations show that, compared to the integer-order model with discrete time-delays, as the fractional order decreases, the critical value for the Hopf bifurcation parameter $\tau=\tau_1+\tau_2$ increases.

\begin{figure}[htbp]
\centering
\begin{tabular}{cc}
\includegraphics[width=0.49\textwidth ]{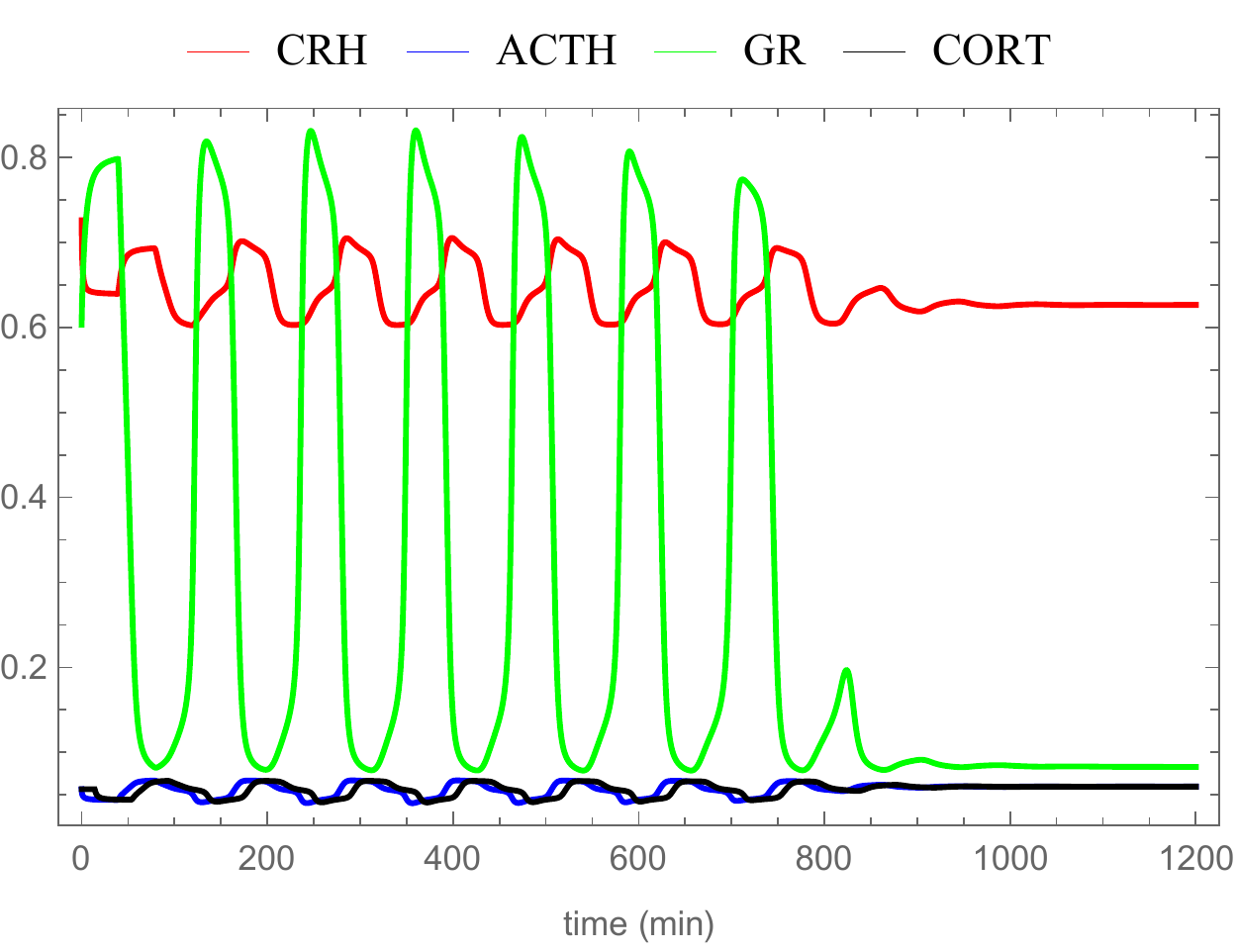} &
\includegraphics[width=0.49\textwidth ]{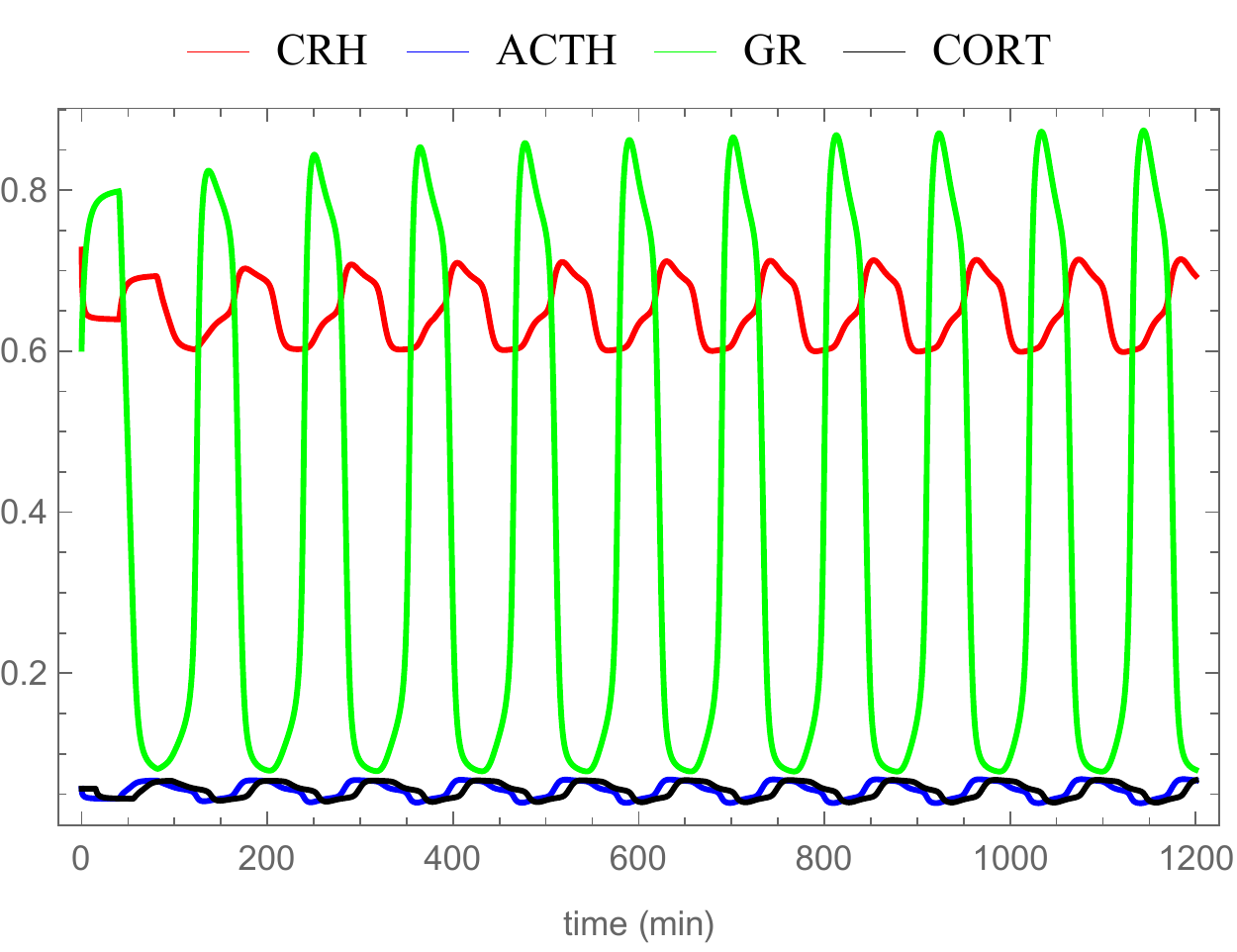}\\
\end{tabular}
\caption{Trajectories of system (\ref{sysHPA4.frac}) with discrete time-delays $\tau_1=25$ min and $\tau_2=14$ min (left) and $\tau_2=15$ min (right), considering an initial condition in a small neighborhood of $E_0$. }
\end{figure}

\section{CONCLUSIONS}
In this paper, we have considered a four-dimensional mathematical model that describes the hypothalamus-pituitary-adrenal axis with the influence of the GR concentration. Due to the fact that the involved processes are not instantaneous, we have incorporated distributed delays, and we have also taken into account fractional-order derivatives for the modelling of memory properties. This approach to the modelling of the biological processes is more realistic because it takes into account the whole past history of the variables, capturing the vital mechanisms of the HPA system. Sufficient conditions for the local asymptotic stability of the equilibrium points have been obtained and a Hopf bifurcation analysis has been undertaken. Numerical simulations reflect the importance of the theoretical results. As a direction for future work, this model can be generalized by considering general feedback functions to account for the interactions within the HPA axis.

\section{ACKNOWLEDGMENTS}
This work was supported by grants of the Romanian National Authority for Scientific Research and Innovation, CNCS-UEFISCDI, project no. PN-II-RU-TE-2014-4-0270 (E. Kaslik) and project no. PN-II-ID-PCE-2011-3-0198 (M. Neam\c{t}u).



%

\end{document}